\documentclass[a4paper]{amsart}
\usepackage{amssymb}
\title[Extremal spectral properties of Lawson surfaces
and Lam\'e equation]{Extremal spectral properties of Lawson tau-surfaces
and the Lam\'e equation}
\author{Alexei V. Penskoi}
\address{Department of Geometry and Topology, 
Faculty of Mathematics and Mechanics, Moscow State University,
Leninskie Gory, GSP-1, 119991, Moscow, Russia\newline \emph{and}
\newline Independent University of Moscow, 
Bolshoy Vlasyevskiy pereulok~11, 119002, Moscow, Russia\newline \emph{and}
\newline Department of Mathematical Modelling (FN-12),
Faculty of Fundamental Sciences,
Bauman Moscow State Technical University,
2-ya~Baumanskaya~5, 105005, Moscow, Russia}
\email{penskoi@mccme.ru}
\date{}
\newtheorem{Proposition}{Proposition}
\newtheorem{Definition}{Definition}
\newtheorem*{Theorem}{Theorem}

\DeclareMathOperator{\Area}{Area}
\DeclareMathOperator{\tr}{tr}
\DeclareMathOperator{\sn}{sn}
\DeclareMathOperator{\cn}{cn}
\DeclareMathOperator{\dn}{dn}
\DeclareMathOperator{\Ec}{Ec}
\DeclareMathOperator{\Es}{Es}
\DeclareMathOperator{\am}{am}
\begin{document}
\begin{abstract}
Extremal spectral properties of Lawson tau-surfaces
are investigated. The Lawson tau-surfaces form a
two-parametric
family of tori or Klein bottles minimally immersed
in the standard unitary three-dimensional sphere. 
A~Lawson tau-surface carries an extremal metric for some
eigenvalue of the Laplace-Beltrami operator.
Using theory of the Lam\'e equation
we find explicitly these extremal
eigenvalues.
\end{abstract}
\maketitle

\section*{Introduction}

Let $M$ be a closed surface and $g$ be
a Riemannian metric on $M.$
Let us consider the associated Laplace-Beltrami
operator $\Delta:C^\infty(M)\longrightarrow C^\infty(M),$
$$
\Delta f=-\frac{1}{\sqrt{|g|}}\frac{\partial}{\partial x^i}%
\left(\sqrt{|g|}g^{ij}\frac{\partial f}{\partial x^j}\right).
$$
It is well-known that the eigenvalues
\begin{equation}\label{eigenvalues}
0=\lambda_0(M,g)<\lambda_1(M,g)\leqslant%
\lambda_2(M,g)\leqslant\lambda_3(M,g)\leqslant\dots
\end{equation}
of $\Delta$ possess the following rescaling property,
$$
\forall t>0\quad\lambda_i(M,tg)=\frac{\lambda_i(M,g)}{t}.
$$
Hence, it is not a good idea to look for
a supremum of the functional $\lambda_i(M,g)$ over the space
of Riemannian metrics $g$ on a fixed surface $M.$
But the functional
$$
\Lambda_i(M,g)=\lambda_i(M,g)\Area(M,g)
$$
is invariant under rescaling transformations $g\mapsto tg.$

It turns out that the question
about the supremum $\sup\Lambda_i(M,g)$
of the functional $\Lambda_i(M,g)$ over the space
of Riemannian metrics $g$ on a fixed surface
$M$ is very difficult and only few results are known.

In 1980 it was proven by Yang and Yau in the
paper~\cite{Yang-Yau} that
for an orientable surface $M$ of genus $\gamma$ the following
inequality holds,
$$
\Lambda_1(M,g)\leqslant 8\pi(\gamma+1).
$$
A generalization of this result for an arbitrary $\Lambda_i$ 
was found in 1993 by Korevaar. It is proven in
the paper~\cite{Korevaar}
that there exists a constant $C$ such that for
any $i>0$ and any compact surface $M$ of genus $\gamma$
the functional $\Lambda_i(M,g)$ is bounded,
$$
\Lambda_i(M,g)\leqslant C(\gamma+1)i.
$$
It should be remarked that in 1994
Colbois and Dodziuk proved in the
paper~\cite{Colbois-Dodziuk}
that for a manifold $M$ of dimension $\dim M\geqslant 3$
the functional $\Lambda_i(M,g)$ 
is not bounded on the space of Riemannian metrics $g$ on $M.$

The functional $\Lambda_i(M,g)$ depends continuously
on the metric $g,$ but this functional
is not differentiable. However, it was shown
in 1973 by Berger in the paper~\cite{Berger} that for analytic
deformations $g_t$ the left and right derivatives
of the functional $\Lambda_i(M,g_t)$ with respect to $t$ exist. 
This led  to the following definition,
see the paper~\cite{Nadirashvili1}
by Nadirashvili (1986) and the paper~\cite{ElSoufi-Ilias1}
by El Soufi and Ilias (2000).

\begin{Definition}
A Riemannian metric $g$ on  a closed surface
$M$ is called extremal metric for the
functional $\Lambda_i(M,g)$ if for any analytic deformation
$g_t$ such that $g_0=g$ the following inequality holds,
$$
\frac{d}{dt}\Lambda_i(M,g_t)%
\left.\vphantom{\raisebox{-0.5em}{.}}\right|_{t=0+}\leqslant0%
\leqslant\frac{d}{dt}\Lambda_i(M,g_t)%
\left.\vphantom{\raisebox{-0.5em}{.}}\right|_{t=0-}.
$$
\end{Definition}

The list of surfaces $M$ and values of index $i$
such that the maximal or at least extremal metrics
for the functional $\Lambda_i(M,g)$ are known is quite short.

\begin{itemize}
\item[$\Lambda_1(\mathbb{S}^2,g)$] Hersch proved
in 1970 in the paper~\cite{Hersch}
that $\sup\Lambda_1(\mathbb{S}^2,g)=8\pi$ and the maximum
is reached on the canonical metric on $\mathbb{S}^2.$
This metric is the unique extremal metric.
\item[$\Lambda_1(\mathbb{R}P^2,g)$] Li and Yau proved
in 1982 in the paper~\cite{Li-Yau}
that $\sup\Lambda_1(\mathbb{R}P^2,g)=12\pi$
and the maximum is reached on the canonical 
metric on $\mathbb{R}P^2.$ This metric is the unique extremal metric.
\item[$\Lambda_1(\mathbb{T}^2,g)$] Nadirashvili
proved in 1996 in the paper~\cite{Nadirashvili1}
that $\sup\Lambda_1(\mathbb{T}^2,g)=\frac{8\pi^2}{\sqrt{3}}$
and the maximum is reached on the flat equilateral torus. El~Soufi
and Ilias proved in 2000 in the paper~\cite{ElSoufi-Ilias1} 
that the only extremal
metric for $\Lambda_1(\mathbb{T}^2,g)$ different from the maximal one
is the metric on the Clifford torus.
\item[$\Lambda_1(\mathbb{K},g)$] Jakobson, Nadirashvili and 
Polterovich proved in 2006 in the paper~\cite{Jakobson-Nadirashvili-Polterovich}
that the metric on a Klein bottle realized
as the Lawson bipolar surface $\tilde{\tau}_{3,1}$
is extremal. El Soufi, Giacomini and 
Jazar proved in the same year in the
paper~\cite{ElSoufi-Giacomini-Jazar}
that this metric is the unique extremal metric and the maximal
one. Here $\sup\Lambda_1(\mathbb{K},g)=12\pi E\left(\frac{2\sqrt{2}}{3}\right),$
where $E$ is a complete elliptic integral of the second kind,
we recall its definition in the end of Introduction.
\item[$\Lambda_2(\mathbb{S}^2,g)$] Nadirashvili
proved in 2002 in the paper~\cite{Nadirashvili2}
that $\sup\Lambda_2(\mathbb{S}^2,g)=16\pi$
and maximum is reached on a singular metric which can be obtained
as the metric on the union of two spheres of equal radius with
canonical metric glued together.
\item[$\Lambda_i(\mathbb{T}^2,g),\Lambda_i(\mathbb{K},g)$]
Let $r,k\in\mathbb{N},$ $0<k<r,$ $(r,k)=1.$
Lapointe studied
bipolar surfaces $\tilde{\tau}_{r,k}$ of Lawson tau-surfaces
$\tau_{r,k}$ and proved the following result published
in 2008 in the paper~\cite{Lapointe}. 
\begin{enumerate}
\item If $rk\equiv 0 \mod 2$ then
$\tilde{\tau}_{r,k}$ is a torus and it carries an extremal metric
for $\Lambda_{4r-2}(\mathbb{T}^2,g).$
\item If $rk\equiv 1 \mod 4$ then
$\tilde{\tau}_{r,k}$ is a torus and it carries an extremal metric
for $\Lambda_{2r-2}(\mathbb{T}^2,g).$
\item If $rk\equiv 3 \mod 4$ then
$\tilde{\tau}_{r,k}$ is a Klein bottle and it carries an extremal metric
for $\Lambda_{r-2}(\mathbb{K},g).$
\end{enumerate}
\end{itemize}

We should also mention the paper~\cite{JLNNP} published in 2005
by Jakobson, Levitin, Nadirashvili, Nigam and Polterovich.
It is shown in this paper using a combination of analytic and numerical tools
that the maximal metric for the first eigenvalue on 
the surface of genus two is the metric on the Bolza surface $\mathcal P$
induced from the canonical metric on the sphere
using the standard covering ${\mathcal P}\longrightarrow\mathbb{S}^2.$ 
In fact, the authors state 
this result as a conjecture, because a part of the argument 
is based on a numerical calculation. 

The goal of the present paper is to study
extremal spectral properties of metrics
on Lawson tau-surfaces $\tau_{m,k}.$

\begin{Definition}\label{Lawson}
A Lawson tau-surface $\tau_{m,k}\subset\mathbb{S}^3$
is defined by the doubly-periodic immersion 
$\Psi_{m,k}:\mathbb{R}^2\longrightarrow\mathbb{S}^3\subset\mathbb{R}^4$
given by the following explicit formula,
\begin{equation}\label{immersion}
\Psi_{m,k}(x,y)=%
(\cos(mx)\cos y, \sin(mx)\cos y, \cos(kx)\sin y,\sin(kx)\sin y).
\end{equation}
\end{Definition}

This family of surfaces is introduced in 1970 by Lawson in
the paper~\cite{Lawson}. He proved that for each unordered pair 
of positive integers $(m,k)$ with $(m,k)=1$ 
the surface $\tau_{m,k}$ is a distinct compact
minimal surface in $\mathbb{S}^3.$
Let us impose the condition $(m,k)=1.$
If both integers
$m$ and $k$ are odd then $\tau_{m,k}$ is a torus.
We call it a Lawson torus. If
one of integers $m$ and $k$ is even then $\tau_{m,k}$ is a
Klein bottle. We call it a Lawson Klein bottle.
Remark that  $m$ and $k$ cannot both
be even due to the identity $(m,k)=1$. 
The torus $\tau_{1,1}$ is the Clifford torus.

Since $\tau_{m,k}\cong\tau_{k,m},$ we can fix a convenient
order of $m$ and $k.$ If $\tau_{m,k}$ is a Lawson torus,
i.e. $m$ and $k$ are both odd, $(m,k)=1,$
let us suppose that always $m>k>0$ except
the special case of the Clifford torus $\tau_{1,1}.$
If $\tau_{m,k}$ is a Lawson Klein bottle,
i.e. one of numbers $m$ and $k$ is even, $(m,k)=1,$
let us suppose that always $m$ is even and
$k$ is odd.

It is clear that the map $\Psi$ has periods
$T_1=(2\pi,0)$ and $T_2=(0,2\pi).$ However,
in the case of a Lawson torus $\tau_{m,k}$
the smallest period lattice is generated by
$T_3=(\pi,\pi)$ and $T_4=(\pi,-\pi).$
Hence, a Lawson torus $\tau_{m,k}$
is isometric to the torus
$\mathbb{R}^2/\{aT_3+bT_4|a,b\in\mathbb{Z}\}$
with the metric induced by the immersion $\Psi.$
We identify these tori.

The torus $\mathbb{R}^2/\{aT_1+bT_2|a,b\in\mathbb{Z}\}$
with the metric induced by the immersion $\Psi$
is a double cover of the Lawson torus $\tau_{m,k}.$
We denote this double cover by $\hat{\tau}_{m,k}.$
When it is necessary
to have uniquely defined coordinates of a point on
$\hat{\tau}_{m,k},$ we consider coordinates 
$(x,y)\in[0,2\pi)\times[-\pi,\pi).$ 
Functions on a Lawson torus $\tau_{m,k}$
are in one-to-one correspondence
with functions on the double cover
$\hat{\tau}_{m,k}$ invariant with respect
to the translation by $T_3.$

In the case of a Lawson Klein bottle $\tau_{m,k}$
the immersion $\Psi$ value is invariant under
transformations
$$
(x,y)\mapsto(x+\pi,-y),\quad(x,y)\mapsto(x,y+2\pi).
$$
When it is necessary
to have uniquely defined coordinates of a point on a Lawson Klein
bottle $\tau_{m,k},$ we consider coordinates 
$(x,y)\in[0,\pi)\times[-\pi,\pi).$

The main result of the present paper
is the following Theorem.
Here $[\alpha]$ denotes the integer part of a real number $\alpha$ and 
$E$ is the complete elliptic integral of the second kind,
we recall its definition in the end of the Introduction.

\begin{Theorem}\label{main}
Let $\tau_{m,k}$ be
a Lawson torus. We can assume that $m,k\equiv 1 \mod 2,$ $(m,k)=1.$
Then the induced metric on its double cover $\hat{\tau}_{m,k}$
is an extremal metric
for the functional $\Lambda_j(\mathbb{T}^2,g),$
where
$$
j=2\left(\left[\sqrt{m^2+k^2}\,\right]+m+k\right)-1.
$$
The corresponding value of the functional is
$$
\Lambda_j(\hat{\tau}_{m,k})=16\pi mE\left(\frac{\sqrt{m^2-k^2}}{m}\right).
$$
The induced metric on $\tau_{m,k}$
is an extremal metric
for the functional $\Lambda_j(\mathbb{T}^2,g),$
where
$$
j=2\left[\frac{\sqrt{m^2+k^2}}{2}\right]+m+k-1.
$$
The corresponding value of the functional is
$$
\Lambda_j(\tau_{m,k})=8\pi mE\left(\frac{\sqrt{m^2-k^2}}{m}\right).
$$

Let $\tau_{m,k}$ be a Lawson Klein bottle. 
We can assume that $m\equiv 0 \mod 2,$
$k\equiv 1 \mod 2,$ $(m,k)=1.$
Then the induced metric on $\tau_{m,k}$
is an extremal metric
for the functional $\Lambda_j(\mathbb{K},g),$
where 
$$
j=2\left[\frac{\sqrt{m^2+k^2}}{2}\right]+m+k-1.
$$
The corresponding value of the functional is
$$
\Lambda_j(\tau_{m,k})=8\pi mE\left(\frac{\sqrt{m^2-k^2}}{m}\right).
$$
\end{Theorem}

We investigate the case of a double cover $\hat{\tau}_{m,k}$ since
in this case we have separation of variables in the
corresponding spectral problem. We reduce the case
of Lawson tori to the case of their double covers.

We already described above the result
by Lapointe from the paper~\cite{Lapointe}.
It implies that
if $rk\equiv 0$ or $rk\equiv 1 \mod 4$ then
the surfaces $\tilde{\tau}_{r,k}$ bipolar to 
Lawson tau-surfaces are tori and the metrics
on these tori are extremal for a
functional $\Lambda_j$ with an {\em even\/} $j.$
In contrast with this result, we prove
that the metrics on the Lawson tori $\tau_{m,k}$
are extremal for a
functional $\Lambda_j$ with an {\em odd\/} $j.$
Hence we provide in this paper extremal metrics
on the torus for eigenvalues such that extremal metrics for them
were not known before. A similar situation with Klein bottles,
Lapointe provided extremal metrics on the Klein bottle
for eigenvalues $\Lambda_j$ with an {\em odd\/} $j$
and we provide extremal metrics on the Klein bottle
for eigenvalues $\Lambda_j$ with an {\em even\/} $j.$
In the same time we should remark that we do not
know at this moment extremal metrics for any $\Lambda_j$
and even an extremal metric for $\Lambda_2(\mathbb{T}^2,g)$ 
is yet unknown.

The proof of the Theorem consists of
several steps. We start by using
a beautiful result relating extremal metrics to minimal
immersions in spheres
proved by El~Soufi and Ilias in the 
paper~\cite{ElSoufi-Ilias2}.
This result reduces calculating $j$
to counting the eigenvalues of the
Laplace-Beltrami operator
$\Delta$ on a Lawson tau-surface $\tau_{m,k}.$
Then we reduce this problem to counting
eigenvalues of a periodic Sturm-Liouville
problem. A crucial idea
of this counting is based on the 
relation to the Lam\'e equation~\eqref{Lame}.

We use on several occasions the complete elliptic
integral of the first kind $K(k)$
and the complete elliptic integral of the second kind $E(k)$
defined by formulae
$$
K(k) = \int\limits_0^1\frac{d\alpha}{\sqrt{1-\alpha^2}\sqrt{1-k^2\alpha^2}},
\quad E(k)= \int\limits_0^1\frac{\sqrt{1-k^2\alpha^2}}{\sqrt{1-\alpha^2}}d\alpha.
$$

\section{Minimal submanifolds of a sphere
and extremal spectral property of their metrics}

Let us recall two important results about
minimal submanifolds of a sphere.
Let $N$ be a $d$-dimensional minimal
submanifold of the sphere $\mathbb{S}^n\subset\mathbb{R}^{n+1}$
of radius~$R.$ Let $\Delta$ be the Laplace-Beltrami operator
on $N$ equipped with the induced metric.

The first result is a classical one. 
Its proof can be found e.g. in the book~\cite{Kobayashi-Nomizu}.

\begin{Proposition}\label{minimal-eigenfunction}
The restrictions $x^1|_N,\dots,x^{n+1}|_N$ on $N$
of the standard coordinate functions of $\mathbb{R}^{n+1}$ 
are eigenfunctions of $\Delta$ with eigenvalue $\frac{d}{R^2}.$
\end{Proposition}

Let us numerate the eigenvalues of $\Delta$ as in formula~\eqref{eigenvalues}
counting them with multiplicities
$$
0=\lambda_0<\lambda_1\leqslant\lambda_2\leqslant\dots\leqslant\lambda_i\leqslant\dots
$$
Proposition~\ref{minimal-eigenfunction} implies that there exists
at least one index $i$ such that $\lambda_i=\frac{d}{R^2}.$
Let $j$ denotes the minimal number $i$ such that $\lambda_i=\frac{d}{R^2}.$
Let us introduce the eigenvalues counting function
$$
N(\lambda)=\#\{\lambda_i|\lambda_i<\lambda\}.
$$
We see that $j=N(\frac{d}{R^2}).$ Remark that we count the
eigenvalues starting from $\lambda_0=0.$

The second result was published in 2008 by 
El~Soufi and Ilias in the paper~\cite{ElSoufi-Ilias2}.
It could be written in the following form.

\begin{Proposition}\label{minimal-extremal}
The metric $g_0$ induced on $N$ by immersion
$N\subset\mathbb{S}^n$ is an extremal metric
for the functional $\Lambda_{N\left(\frac{d}{R^2}\right)}(N,g).$
\end{Proposition}

We should remark that isometric immersions by eigenfunctions
of the Laplace-Beltrami operator were studied
in the paper~\cite{Takahashi} by Takahashi. The results
of Takahashi are used by El~Soufi and Ilias in the 
paper~\cite{ElSoufi-Ilias2}.

Proposition~\ref{minimal-extremal} implies immediately the following one.

\begin{Proposition}\label{reduction}
The  metric $g_0$ induced on a Lawson 
torus or Klein bottle 
$\tau_{m,k}$ by its immersion $\tau_{m,k}\subset\mathbb{S}^3$
is an extremal metric for the functional 
$\Lambda_{N(2)}(\mathbb{T}^2,g)$ or $\Lambda_{N(2)}(\mathbb{K},g),$
respectively. The similar statement holds also for the double
cover~$\hat{\tau}_{m,k}.$
\end{Proposition}

\noindent{\bf Proof.} As we mentioned
in the Introduction, Lawson proved in the paper~\cite{Lawson}
that $\tau_{m,k}$ is a complete minimal 
surface in the sphere $\mathbb{S}^3$
of radius $1.$
The statement follows immediately from  
Proposition~\ref{minimal-extremal}, where $R=1$ and $d=2.$
$\Box$

Proposition~\ref{reduction} is crucial for this paper. It reduces
investigation of extremal spectral properties of the Lawson
tau-surfaces to counting eigenvalues $\lambda_i$ of the Laplace-Beltrami
operator such that $\lambda_i<2.$

\section{Eigenvalues of the Laplace-Beltrami
operator on Lawson tau-surfaces and auxiliary periodic
Sturm-Liouville problem}

\begin{Proposition}
Let $\tau_{m,k}\subset\mathbb{S}^3$ be a Lawson tau-surface
and 
$$
p(y)=\sqrt{k^2+(m^2-k^2)\cos^2y}.
$$
Then the induced metric is equal to
$$
p^2(y)dx^2+dy^2
$$
and the Laplace-Beltrami operator is given by 
the following formula,
$$
\Delta f=-\frac{1}{p^2(y)}\frac{\partial^2 f}{\partial x^2}-%
\frac{1}{p(y)}\frac{\partial}{\partial y}%
\left(p(y)\frac{\partial f}{\partial y}\right).
$$
The same holds for a double cover $\hat{\tau}_{m,k}.$
\end{Proposition}

\noindent{\bf Proof} is by direct calculation. $\Box$

Counting eigenvalues is known to be a difficult problem.
Fortunately, we can reduce this problem to a one-dimensional
one. Let us recall that a function $\varphi(y)$
is called $\pi$-antiperiodic if $\varphi(y+\pi)=-\varphi(y).$

\begin{Proposition}\label{reduction-proposition}
Let $\varphi(l,y)$ be a solution of
a periodic Sturm-Liouville problem
\begin{gather}
-\frac{1}{p(y)}\frac{d}{dy}\left(p(y)\frac{d\varphi(y)}{dy}\right)+%
\left(\frac{l^2}{p^2(y)}-\lambda\right)\varphi(y)=0,\label{Sturm-Liouville}\\
\varphi(y+2\pi)\equiv\varphi(y).\label{periodic-phi}
\end{gather}
Let $\hat{\tau}_{m,k}$ be a double cover of a Lawson torus.
Then functions
\begin{equation}\label{phi-cos}
\varphi(l,y)\cos(lx),\quad l=0,1,2,\dots,
\end{equation}
and
\begin{equation}\label{phi-sin}
\varphi(l,y)\sin(lx),\quad l=1,2,3,\dots,
\end{equation}
form a basis in the space of eigenfunctions of the Laplace-Beltrami
operator $\Delta$ with eigenvalue $\lambda.$

Let $\tau_{m,k}$ be a Lawson torus.
Then functions
$$
\varphi(l,y)\cos(lx),\quad l=0,2,4\dots,
$$
and
$$
\varphi(l,y)\sin(lx),\quad l=2,4,6\dots,
$$
where $\varphi(l,y)$ is a $\pi$-periodic solution of
the Sturm-Liouville
problem~\eqref{Sturm-Liouville},
and functions
$$
\varphi(l,y)\cos(lx),\quad l=1,3,5\dots,
$$
and
$$
\varphi(l,y)\sin(lx),\quad l=1,3,5\dots,
$$
where $\varphi(l,y)$ is a $\pi$-antiperiodic solution of
the Sturm-Liouville
problem~\eqref{Sturm-Liouville},
form a basis in the space of eigenfunctions of the Laplace-Beltrami
operator $\Delta$ with eigenvalue $\lambda.$

Let $\tau_{m,k}$ be a Lawson Klein bottle.
Then functions
$$
\varphi(l,y)\cos(lx),\quad l=0,2,4\dots,
$$
and
$$
\varphi(l,y)\sin(lx),\quad l=2,4,6\dots,
$$
where $\varphi(l,y)$ is an even solution of
the periodic Sturm-Liouville
problem~\eqref{Sturm-Liouville}, \eqref{periodic-phi},
and functions
$$
\varphi(l,y)\cos(lx),\quad l=1,3,5\dots,
$$
and
$$
\varphi(l,y)\sin(lx),\quad l=1,3,5\dots,
$$
where $\varphi(l,y)$ is an odd solution of
the periodic Sturm-Liouville
problem~\eqref{Sturm-Liouville}, \eqref{periodic-phi},
form a basis in the space of eigenfunctions of the Laplace-Beltrami
operator $\Delta$ with eigenvalue $\lambda.$
\end{Proposition}

\noindent{\bf Proof.}
Let us remark that $\Delta$ commutes 
with $\frac{\partial}{\partial x}.$ It follows
that $\Delta$ has a basis
of eigenfunctions of the form $\varphi(l,y)\cos(lx)$
and $\varphi(l,y)\sin(lx)$. Substituting 
these eigenfunctions
into the formula $\Delta f=\lambda f,$ we obtain
equation~\eqref{Sturm-Liouville}. 
However, these solutions should be invariant under
transformations
$$
(x,y)\mapsto(x+2\pi,y),\quad(x,y)\mapsto(x,y+2\pi)
$$
in the case of a double cover $\hat{\tau}_{m,k},$
the same transformations plus
$$
(x,y)\mapsto(x+\pi,y+\pi)
$$
in the case of a Lawson torus,
and
$$
(x,y)\mapsto(x+\pi,-y),\quad(x,y)\mapsto(x,y+2\pi)
$$
in the case of a Lawson Klein bottle.

These conditions
imply the periodicity conditions~\eqref{periodic-phi}
and the conditions on $l$ and parity or (anti-)periodicity
of $\varphi.$
$\Box$

We consider $l$ as a parameter in equation~\eqref{Sturm-Liouville}
and in its solutions $\varphi(l,y).$ For example, we assume
$l$ fixed and consider $y$ as an independent 
variable when we discuss
zeroes of the function $\varphi(l,y).$

Let us rewrite equation~\eqref{Sturm-Liouville}
in the standard form of a Sturm-Liouville problem,
\begin{equation}\label{SL-standard}
(p(y)\varphi'(y))'+\left(\lambda p(y)-\frac{l^2}{p(y)}\right)\varphi(y)=0.
\end{equation}
Since $p(y)>0,$ the following classical result holds, 
see e.g. the book~\cite{Coddington-Levinson}.

\begin{Proposition}\label{SL-properties}
There are an infinite number of eigenvalues $\lambda_i(l)$ 
of the periodic Sturm-Liouville problem~\eqref{Sturm-Liouville},
\eqref{periodic-phi}.
The eigenvalues $\lambda_i(l)$  form a sequence such that
$$
\lambda_0(l)<\lambda_1(l)\leqslant\lambda_2(l)<%
\lambda_3(l)\leqslant\lambda_4(l)<%
\lambda_5(l)\leqslant\lambda_6(l)<\dots
$$
For $\lambda=\lambda_0(l)$ there exists a unique
(up to
a multiplication by a non-zero constant)
eigenfunction $\varphi_0(l,y).$ If 
$\lambda_{2i+1}(l)<\lambda_{2i+2}(l)$
for some $i\geqslant0$ then there is a unique (up to
a multiplication by a non-zero constant) eigenfunction 
$\varphi_{2i+1}(l,y)$ with eigenvalue $\lambda=\lambda_{2i+1}(l)$
of multiplicity one
and there is a unique (up to
a multiplication by a non-zero constant) eigenfunction 
$\varphi_{2i+2}(l,y)$ with eigenvalue $\lambda=\lambda_{2i+2}(l)$
of multiplicity one.
If  $\lambda_{2i+1}(l)=\lambda_{2i+2}(l)$
then there are two independent eigenfunctions
$\varphi_{2i+1}(l,y)$ and $\varphi_{2i+2}(l,y)$ 
with eigenvalue $\lambda=\lambda_{2i+1}(l)=\lambda_{2i+1}(l)$
of multiplicity two.

The eigenfunction $\varphi_0(l,y)$ has no zeroes on $[0,2\pi).$
The eigenfunctions 
$\varphi_{2i+1}(l,y)$ and $\varphi_{2i+2}(l,y),$
$i\geqslant0,$ each  
have exactly $2i+2$ zeroes on $[0,2\pi).$
\end{Proposition}

As usual, we can rewrite our periodic
Sturm-Liouville problem as a periodic
problem for a Hill equation.

\begin{Proposition}\label{SL-Hill}
Equation~\eqref{Sturm-Liouville}
is equivalent to a Hill equation
\begin{equation}\label{Hill}
-z''(y)+V(l,y)z(y)=\lambda z(y),
\end{equation}
where
$$
V(l,y)=\frac{l^2}{p(y)^2}+\frac{1}{4}\left(\frac{p'(y)}{p(y)}\right)^2+%
\frac{1}{2}\left(\frac{p'(y)}{p(y)}\right)'.
$$
Periodic boundary condition~\eqref{periodic-phi}
for equation~\eqref{Sturm-Liouville}
is equivalent to the periodic
boundary condition
\begin{equation}\label{periodic-z}
z(y+2\pi)\equiv z(y).
\end{equation}
\end{Proposition}

\noindent{\bf Proof.} A direct calculation
shows that the change of
variable
$z(y)=\sqrt{p(y)}\varphi(y)$
transforms equation~\eqref{Sturm-Liouville}
into equation~\eqref{Hill}. Since the
function $p(y)$
is $2\pi$-periodic, boundary conditions~\eqref{periodic-phi}
and~\eqref{periodic-z} are equivalent. $\Box$

Since $p(y)$ is an even $\pi$-periodic function,
the following propositions hold.

\begin{Proposition}\label{SL-parity}
The solution $\varphi_0(l,y)$ is even.
If $\lambda_{i}(l),$ $i>0,$ is of multiplicity one
then the solution $\varphi_{i}(l,y)$
is even or odd.
If  $\lambda_{2i+1}(l)=\lambda_{2i+2}(l)$
then two independent eigenfunctions
$\varphi_{2i+1}(l,y)$ and $\varphi_{2i+2}(l,y)$ 
with eigenvalue $\lambda_{2i+1}(l)=\lambda_{2i+1}(l)$
of multiplicity two
could be chosen in such a way that one of them is even
and another is odd.
\end{Proposition}

\noindent{\bf Proof.} Applying Theorem 1.1
from the book~\cite{Magnus-Winkler} to
equation~\eqref{Hill} with periodic
condition~\eqref{periodic-z}
and returning then back from $z(y)$
to $\varphi(y),$ we obtain the statement.
$\Box$

\begin{Proposition}\label{SL-pi}
The solution $\varphi_0(l,y)$ is $\pi$-periodic.
The solutions $\varphi_{2i+1}(l,y)$ and 
$\varphi_{2i+2}(l,y)$ are $\pi$-periodic
if $i$ is odd and $\pi$-antiperiodic
if $i$ is even.
\end{Proposition}

\noindent{\bf Proof.} Applying Theorem 3.1
from Chapter VIII of the book~\cite{Coddington-Levinson}
with period $\pi$ and $2\pi$ and comparing,
we immediately obtain the statement.
$\Box$

It is easy now to establish a relation between
the multiplicities of the eigenvalues of the operator $\Delta$
and the eigenvalues $\lambda_i(l)$
of the periodic Sturm-Liouville 
problem~\eqref{Sturm-Liouville},~\eqref{periodic-phi}.
This relation permits us to express the quantity $N(2)$
in terms of the eigenvalues $\lambda_i(l).$

\begin{Proposition}\label{how-to-count}
In the case of a double cover $\hat{\tau}_{m,k}$ 
of a Lawson torus we have  
\begin{equation}\label{N2cover}
N(2)=\#\{\lambda_i(0)|\lambda_i(0)<2\}+%
2\#\{\lambda_i(l)|\lambda_i(l)<2,l>0,l\in\mathbb{Z}\}.
\end{equation}
In the case of a Lawson torus $\tau_{m,k}$
we have 
\begin{gather}
N(2)=1+%
2\#\{\lambda_0(l)|\lambda_0(l)<2,l>0,l\in\mathbb{Z},%
l\mbox{\ is even}\}+\label{N2torus}\\
+\#\{\lambda_{2i+1}(0)|\lambda_{2i+1}(0)<2,i\mbox{\ is odd}\}%
+\#\{\lambda_{2i+2}(0)|\lambda_{2i+2}(0)<2,i\mbox{\ is odd}\}+\notag\\
+2\#\{\lambda_{2i+1}(l)|\lambda_{2i+1}(l)<2,l>0,l\in\mathbb{Z},%
l\mbox{\ is even},i\mbox{\ is odd}\}+\notag\\
+2\#\{\lambda_{2i+2}(l)|\lambda_{2i+2}(l)<2,l>0,l\in\mathbb{Z},%
l\mbox{\ is even},i\mbox{\ is odd}\}+\notag\\
+2\#\{\lambda_{2i+1}(l)|\lambda_{2i+1}(l)<2,l>0,l\in\mathbb{Z},%
l\mbox{\ is odd},i\mbox{\ is even}\}+\notag\\
+2\#\{\lambda_{2i+2}(l)|\lambda_{2i+2}(l)<2,l>0,l\in\mathbb{Z},%
l\mbox{\ is odd},i\mbox{\ is even}\}.\notag
\end{gather}
In the case of a Lawson Klein bottle $\tau_{m,k}$
we have 
\begin{gather}
N(2)=\#\{\lambda_i(0)|\lambda_i(0)<2,\varphi_i(0,y)\mbox{\ is even}\}+\label{N2Klein}\\
+2\#\{\lambda_i(l)|\lambda_i(l)<2,l>0,l\in\mathbb{Z},l\mbox{\ is even},%
\varphi_i(l,y)\mbox{\ is even}\}+\notag\\
+2\#\{\lambda_i(l)|\lambda_i(l)<2,l>0,l\in\mathbb{Z},l\mbox{\ is odd},%
\varphi_i(l,y)\mbox{\ is odd}\}.\notag
\end{gather}
\end{Proposition}

\noindent{\bf Proof.} Let $\tau_{m,k}$ be a Lawson torus.
The eigenvalue $\lambda_i(0)$ gives
exactly one basis eigenfunction 
$$
\varphi_i(0,y)\cos(0x)=\varphi_i(0,y)
$$
of the operator~$\Delta.$
It follows that
each $\lambda_i(0)$ corresponds to one eigenvalue
$\lambda_q=\lambda_i(0)$ of the Laplace-Beltrami operator $\Delta.$
The eigenvalue $\lambda_i(l),$ $l>0,$ gives
exactly two basis eigenfunctions 
$$\varphi_i(l,y)\cos(lx)\quad
\mbox{and}\quad \varphi_i(l,y)\sin(lx)
$$ of the Laplace-Beltrami operator~$\Delta.$
It follows that 
for $l>0$ each $\lambda_i(l)$ corresponds to two eigenvalues
$\lambda_q=\lambda_{q+1}=\lambda_i(l)$ of the Laplace-Beltrami operator $\Delta.$
This implies formula~\eqref{N2cover}. 

The cases of a Lawson torus or a Lawson Klein bottle
are similar, but we should take into account parity or
(anti-)periodicity of $\varphi_i(l,y)$ 
see Propositions~\ref{reduction-proposition},
\ref{SL-parity} and~\ref{SL-pi}.
The term $1$ in formula~\eqref{N2torus} counts
$\lambda_0(0).$
$\Box$

Propositions~\ref{reduction-proposition} 
and~\ref{how-to-count} permit
us to investigate the simplest case of the
Lawson torus $\tau_{1,1},$ i.e. the Clifford
torus.

\begin{Proposition}\label{1-1}
The metric induced on the Lawson
torus $\tau_{1,1}$ (i.e. the Clifford
torus) is an extremal metric for the
functional $\Lambda_1(\mathbb{T}^2,g).$
The corresponding value of this
functional is $\Lambda_1(\tau_{1,1})=4\pi^2=39.479\dots$ 

The metric induced on the double cover
$\hat{\tau}_{1,1}$ of the Clifford
torus $\tau_{1,1}$ is an extremal metric for the
functional $\Lambda_5(\mathbb{T}^2,g).$
The corresponding value of this
functional is $\Lambda_5(\hat{\tau}_{1,1})=8\pi^2=78.957\dots$ 

\end{Proposition}

\noindent{\bf Proof.} We have $m=k=1,$ hence
$p(y)\equiv1.$ Equation~\eqref{Sturm-Liouville}
becomes the equation
\begin{equation}\label{SL-1-1}
-(\varphi(y))''+(l^2-\lambda)\varphi(y)=0.
\end{equation}
It is well-known that the eigenvalues $\lambda_i(l)$
of this equation with the periodic boundary
conditions~\eqref{periodic-phi}
satisfy the relation 
$$
\lambda_i(l)-l^2=n^2,\quad n\in\mathbb{Z}.
$$
Hence we have the following eigenvalues
of the periodic Sturm-Liouville problem
that are less then $2,$
\begin{equation*}
\begin{split}
&\lambda_0(0)=0,\quad\lambda_1(0)=1,\quad\lambda_2(0)=1,\\
&\lambda_0(1)=1.
\end{split}
\end{equation*}
It follows from formula~\eqref{N2torus} that 
$N(2)=1$ and the induced metric is
extremal for the functional 
$\Lambda_{1}(\mathbb{T}^2,g).$
Since $\lambda_1=2$ and $\Area(\tau_{1,1})=2\pi^2,$
the value of this functional is $\Lambda_1(\tau_{1,1})=4\pi^2.$

The case of the double cover $\hat{\tau}_{1,1}$ is similar.
$\Box$

This result is well-known and the example 
of $\tau_{1,1}$ is a "toy example" because
in this simplest case the corresponding equation~\eqref{SL-1-1}
is exactly solvable. It is not the case for other Lawson tau-surfaces.
However, we can find $N(2)$ by investigating
the structure of the eigenvalues $\lambda_i(l)$ of
the auxiliary Sturm-Liouville 
problems~\eqref{Sturm-Liouville},~\eqref{periodic-phi}.

\section{Magnus-Winkler-Ince equation and eigenvalues of multiplicity $2$}

It turns out that it is important
to investigate eigenvalues of multiplicity $2$ 
of the periodic Sturm-Liouville 
problem~\eqref{Sturm-Liouville}, \eqref{periodic-phi}.
The problem of existence
of such eigenvalues is usually called a coexistence
problem since this means that for such an eigenvalue
two independent $2\pi$-periodic solutions exist.

Since the Lawson torus $\tau_{1,1}$ is already investigated,
we assume that $m\ne k$ in this section.

Let us remark that our equation~\eqref{Sturm-Liouville}
can be written as
\begin{equation}\label{MWI}
(1+a\cos(2y))\varphi''(y)+b\sin(2y)\varphi'(y)+(c+d\cos(2y))\varphi(y)=0,
\end{equation}
where
\begin{equation}\label{MWI-values}
a=\frac{m^2-k^2}{m^2+k^2},\quad b=-\frac{m^2-k^2}{m^2+k^2},\quad
c=\lambda_i(l)-\frac{2l^2}{m^2+k^2},\quad 
d=\lambda_i(l)\frac{m^2-k^2}{m^2+k^2}.
\end{equation}
Equation~\eqref{MWI} is called a 
Magnus-Winkler-Ince (MWI) equation. Theory
of this equation is interesting for us
because the coexistence problem for the MWI
equation was intensively studied, see
the book~\cite{Magnus-Winkler},
and solved completely in some terms by Volkmer
in 2003, see 
the paper~\cite{Volkmer}.
We need the following 
result~\cite[Theorem~7.1]{Magnus-Winkler}.

\begin{Proposition}\label{MWI-statement}
If the MWI equation~\eqref{MWI} has two linearly
independent solutions of period $2\pi$ then
the polynomial
\begin{equation}\label{Q-polynomial}
Q^*(\mu)=a(2\mu-1)^2-b(2\mu-1)-d
\end{equation}
vanishes for one of the values
of $\mu=0,\pm1,\pm2,\dots.$
\end{Proposition}

This proposition implies immediately  the following Proposition. 

\begin{Proposition}\label{lambda-possible} Let $m\ne k.$
The only possible eigenvalue $\lambda_i(l)$ such that
$\lambda_i(l)<6$ and $\lambda_i(l)$ has multiplicity
$2$ is $\lambda_i(l)=2.$
\end{Proposition}

\noindent{\bf Proof.}
If $\lambda_i(l)$ is an eigenvalue of multiplicity $2$
then by Proposition~\ref{MWI-statement}
the polynomial $Q^*(\mu)$ has integer root~$\mu_0.$ 
Substituting formulae~\eqref{MWI-values} 
into equation~\eqref{Q-polynomial},
we obtain
$$
Q^*(\mu)=\frac{m^2-k^2}{m^2+k^2}\left[(2\mu-1)^2+(2\mu-1)-\lambda_i(l)\right].
$$
It follows that $Q^*(\mu)=0$ is equivalent to
$$
\lambda_i(l)=2\mu(2\mu-1).
$$
Then 
$\lambda_i(l)=2\mu_0(2\mu_0-1),$ $\mu_0\in\mathbb{Z}.$
The only possible 
eigenvalues $0\leqslant\lambda_i(l)<6$
of this form are $\lambda_i(l)=0$ and $\lambda_i(l)=2.$
We know that $\lambda=0$ is the eigenvalue of multiplicity~$1$ 
of $\Delta$ because there is no harmonic functions
on a compact surface except constants. Hence, $\lambda_i(l)=0$ is
excluded and $\lambda_i(l)=2$ is the only possibility.
$\Box$

We know from Proposition~\ref{minimal-eigenfunction}
that the restrictions on $\tau_{m,k}$
of the coordinate functions $x^1,\dots,x^4$
are eigenfunctions of $\Delta$
with eigenvalue $2.$ The explicit formula of
the immersion~\eqref{immersion} gives these functions,
\begin{equation}\label{eigenfunctions-immersion}
\cos(mx)\cos y, \quad \sin(mx)\cos y, \quad 
\cos(kx)\sin y,\quad \sin(kx)\sin y.
\end{equation}

\begin{Proposition}\label{lambda-2-1} Let $m\ne k.$

The function $\cos y$ is an eigenfunction of the 
periodic Sturm-Liouville
problem~\eqref{Sturm-Liouville}, \eqref{periodic-phi} with $l=m$
and the corresponding eigenvalue is equal to $2.$ 
This eigenvalue has multiplicity~$1.$

The function $\sin y$ is an eigenfunction of the 
periodic Sturm-Liouville
problem~\eqref{Sturm-Liouville}, \eqref{periodic-phi} with $l=k$
and the corresponding eigenvalue is equal to $2.$
This eigenvalue has multiplicity~$1.$
\end{Proposition}

\noindent{\bf Proof.}
The eigenfunctions~\eqref{eigenfunctions-immersion}
are of the form~\eqref{phi-cos}
and~\eqref{phi-sin}. This implies that $\cos y$
is an eigenfunction of the periodic Sturm-Liouville
problem~\eqref{Sturm-Liouville}, \eqref{periodic-phi} with $l=m$
and the corresponding eigenvalue is equal to $2.$ 
In the same way, $\sin y$ is an eigenfunction 
of the periodic Sturm-Liouville
problem~\eqref{Sturm-Liouville}, \eqref{periodic-phi} 
with $l=k$ and the corresponding eigenvalue is equal to $2.$ 

Let us consider the case
$l=m$ and $\varphi(y)=\cos y.$
The standard  argument with the Wronskian shows
that a linearly independent with $\cos y$ solution
of equation~\eqref{Sturm-Liouville}
can locally be presented in the form
$$
\cos y\int\limits_{y_0}^y\frac{d\xi}{\cos^2\xi\sqrt{k^2+(m^2-k^2)\cos^2\xi}}.
$$
This integral has singularities at 
$y=\pm\frac{\pi}{2},\pm\frac{3\pi}{2},\dots$
and it is not possible to present a second solution in this
form globally.
Let us define a function
$F:(-\frac{\pi}{2},\frac{\pi}{2})\longrightarrow\mathbb{R}$
by the formula
$$
F(y)=\cos y\int\limits_0^y\frac{d\xi}{\cos^2\xi\sqrt{k^2+(m^2-k^2)\cos^2\xi}}
$$
and a function $G:(\frac{\pi}{2},\frac{3\pi}{2})\longrightarrow\mathbb{R}$
by the formula
$$
G(y)=\cos y\int\limits_\pi^y\frac{d\xi}{\cos^2\xi\sqrt{k^2+(m^2-k^2)\cos^2\xi}}.
$$
It is easy to see that the function $F$ is an odd function.
Let us remark that
$$
\lim_{y\to\frac{\pi}{2}-}F(y)=\frac{1}{k},
\lim_{y\to\frac{\pi}{2}-}F'(y)=%
\frac{1}{k}\left(E\left(\frac{\sqrt{k^2-m^2}}{k}\right)-%
K\left(\frac{\sqrt{k^2-m^2}}{k}\right)\right)\ne0.
$$
The expression for $\lim\limits_{y\to\frac{\pi}{2}-}F'(y)$
is long, let us denote it by $Z$ in order to shorten the notation.

The identity $G(y)=F(\pi-y)$
implies 
$$
\lim_{y\to\frac{\pi}{2}+}G(y)=\lim_{y\to\frac{\pi}{2}-}F(y)=\frac{1}{k},\quad
\lim_{y\to\frac{\pi}{2}+}G'(y)=-\lim_{y\to\frac{\pi}{2}-}F'(y)=-Z\ne0.
$$

It follows that we can define a
linearly independent with $\cos y$ solution
on $[-\frac{\pi}{2},\frac{3\pi}{2}]$
by the formula
$$
\psi(y)=\begin{cases}
\frac{1}{k},&\text{if $x=-\frac{\pi}{2};$}\\
F(y),&\text{if $x\in(-\frac{\pi}{2},\frac{\pi}{2});$}\\
\frac{1}{k},&\text{if $x=\frac{\pi}{2};$}\\
-2Z\cos y%
+G(y),&\text{if $x\in(\frac{\pi}{2},\frac{3\pi}{2});$}\\
\frac{1}{k},&\text{if $x=\frac{3\pi}{2}.$}\\
\end{cases}
$$
This is a smooth solution. However, $\psi(y)$
is not periodic because
$$
\lim\limits_{y\to-\frac{\pi}{2}+}\psi'(y)=%
\lim\limits_{y\to-\frac{\pi}{2}+}F'(y)=%
\lim\limits_{y\to\frac{\pi}{2}-}F'(y)=Z\ne0
$$
and
$$
\lim\limits_{y\to\frac{3\pi}{2}-}\psi'(y)=%
-2Z+%
\lim\limits_{y\to\frac{3\pi}{2}-}G'(y)=%
-3Z\ne0.
$$
This implies that the eigenvalue $\lambda_i(m)=2$
corresponding to the eigenfunction $\cos y$
has multiplicity $1.$

The case $l=k$ and $\sin y$ can be investigated in the same way.
$\Box$

\section{Properties of eigenvalue $\lambda_i(l)$ as a function of $l$}

Let us now consider the periodic Sturm-Liouville
problem~\eqref{Sturm-Liouville}, \eqref{periodic-phi}
not only for integer values of $l$ but also for real 
values of~$l.$

Let us recall a little bit of the general theory
of a Sturm-Liouville problem,
see e.g. the textbook~\cite{Coddington-Levinson}.
Let us rewrite equation~\eqref{SL-standard}
as
\begin{equation}\label{Sturm-Liouville-1}
\left(\frac{p(y)}{m}\varphi'(y)\right)'+%
\left(\lambda\frac{p(y)}{m}-\frac{l^2}{mp(y)}\right)\varphi(y)=0.
\end{equation}
This form is more standard because $\frac{p(0)}{m}=1.$
Let $\Phi(\lambda,y)$ and $\Psi(\lambda,y)$ denote two solutions of
equation~\eqref{Sturm-Liouville-1} such that
$$
\Phi(\lambda,0)=1,\quad\Phi'(\lambda,0)=0,
\quad\Psi(\lambda,0)=0,\quad\Psi'(\lambda,0)=1.
$$
They form a basis in the space of solutions of 
equation~\eqref{Sturm-Liouville-1}. The matrix of the shift
operator $(T\varphi)(y)=\varphi(y+2\pi)$ in this basis
is equal to
$$
\hat{T}(\lambda)=\left(\begin{array}{cc}
\Phi(\lambda,2\pi)&\Psi(\lambda,2\pi)\\
\Phi'(\lambda,2\pi)&\Psi'(\lambda,2\pi)
\end{array}\right).
$$
The conservation law for the Wronskian implies 
$\det\hat{T}(\lambda)=1.$
Then the eigenvalues $\mu$ of the matrix $\hat{T}(\lambda)$ 
are roots of the polynomial
\begin{equation}\label{T-char}
\mu^2-\tr\hat{T}(\lambda)\mu+\det\hat{T}(\lambda)%
=\mu^2-\tr\hat{T}(\lambda)\mu+1.
\end{equation}
It is clear that $\lambda$ is an eigenvalue of
the periodic Sturm-Liouville problem
for equation~\eqref{Sturm-Liouville-1}
if and only
if $\mu=1$ is a root of polynomial~\eqref{T-char}.
This is equivalent to the equation
$\tr\hat{T}(\lambda)=2.$

Let us denote $\tr\hat{T}(\lambda)$ by $f(l,\lambda),$ i.e.
$$
f(l,\lambda)=\Phi(\lambda,2\pi)+\Psi'(\lambda,2\pi).
$$
Then $\lambda_i(l)$ is defined implicitly by the equation
$$
f(l,\lambda)=2.
$$
It is known (see e.g. the textbook~\cite{Coddington-Levinson})
that if $\lambda_i(l)$ is an eigenvalue of multiplicity~$1$ 
then $\frac{\partial f}{\partial\lambda}(l,\lambda_i(l))\ne0.$

\begin{Proposition}\label{lambda-l}
Let us fix $i.$ If $\lambda_i(l)$ has multiplicity $1$
for all $l\in(0,l_1)$ then
$\lambda_i(l)$ is a strictly increasing function
on $(0,l_1)$.
\end{Proposition}

\noindent{\bf Proof.}
Let us introduce a new parameter $\varkappa=l^2.$
We use $\varkappa$ since equation~\eqref{Sturm-Liouville-1}
depends on $\varkappa$ in a linear way.
If $\lambda_i(l)$ has multiplicity $1$ 
then $\frac{\partial f}{\partial\lambda}(l,\lambda_i(l))\ne0.$
This implies that
$\frac{\partial f}{\partial\varkappa}%
(\varkappa,\lambda_i(\varkappa))\ne0$
because 
$\frac{\partial\varkappa}{\partial l}=%
\frac{\partial l^2}{\partial l}=2l\ne0.$
It follows from the implicit function theorem that
$$
\frac{\partial\lambda_i(\varkappa)}{\partial\varkappa}=%
-\frac{\frac{\partial f}{\partial\lambda}(\varkappa,\lambda_i(\varkappa))}%
{\frac{\partial f}{\partial\varkappa}(\varkappa,\lambda_i(\varkappa))}.
$$
It is known (see e.g. the textbook~\cite{Coddington-Levinson})
that
$$
\frac{\partial f}{\partial\lambda}(\varkappa,\lambda_i(\varkappa))=%
\int\limits_0^{2\pi}\left[\Psi^2(\tau)\Phi'(2\pi)+\Psi(\tau)\Phi(\tau)%
(\Phi(2\pi)-\Psi'(2\pi))-\Phi^2(\tau)\Psi(2\pi)\right]\frac{p(\tau)}{m}d\tau
$$
and $\frac{\partial f}{\partial\lambda}(\varkappa,\lambda_i(\varkappa))>0$
for odd $i$ and $\frac{\partial f}{\partial\lambda}(\varkappa,\lambda_i(\varkappa))<0$
for even $i.$

One can prove in a similar way that
$$
\frac{\partial f}{\partial\varkappa}(\varkappa,\lambda_i(\varkappa))=%
-\int\limits_0^{2\pi}\left[\Psi^2(\tau)\Phi'(2\pi)+\Psi(\tau)\Phi(\tau)%
(\Phi(2\pi)-\Psi'(2\pi))-\Phi^2(\tau)\Psi(2\pi)\right]\frac{d\tau}{mp(\tau)}
$$
and $\frac{\partial f}{\partial\varkappa}(\varkappa,\lambda_i(\varkappa))<0$
for odd $i$ and 
$\frac{\partial f}{\partial\varkappa}(\varkappa,\lambda_i(\varkappa))>0$
for even $i.$

It follows that 
$$
\frac{\partial\lambda_i(\varkappa)}{\partial\varkappa}=%
-\frac{\frac{\partial f}{\partial\lambda}(\varkappa,\lambda_i(\varkappa))}%
{\frac{\partial f}{\partial\varkappa}(\varkappa,\lambda_i(\varkappa))}>0
$$
for any parity of $i.$ $\Box$

\section{The Lam\'e equation}\label{Lame-section}

In this Section we recall some properties of
the Lam\'e equation usually written as
\begin{equation}\label{Lame}
\frac{d^2\varphi}{dz^2}+(h-n(n+1)[\hat{k}\sn(z)]^2)\varphi=0.
\end{equation}
see e.g. the book~\cite{BatemanErdelyi}
or the book~\cite{Arscott}. We denote
the modulus of the elliptic function $\sn z$ by $\hat{k}$ since
we already use a letter $k$ in $\tau_{m,k}.$

The Lam\'e equation could be written in different forms,
we will use the trigonometric form of the Lam\'e
equation
\begin{equation}\label{Lame-trig}
[1-(\hat{k}\sin y)^2]\frac{d^2\varphi}{dy^2}-%
\hat{k}^2\sin y\cos y\frac{d\varphi}{dy}+%
[h-n(n+1)(\hat{k}\sin y)^2]\varphi=0.
\end{equation}
Equation~\eqref{Lame-trig} could be obtained
from equation~\eqref{Lame} using the change of variable
\begin{equation}\label{sn-sin}
\sn z=\sin y\quad\Longleftrightarrow\quad y=\am z,
\end{equation}
where $\am z$ is Jacobi amplitude function,
see e.g. the book~\cite[Section~13.9]{BatemanErdelyi}.
This trigonometric form of the Lam\'e equation
is used in the book~\cite{Arscott}.
The change of variable $\sn z=\cos y$
leads to another trigonometric form~\eqref{Lame-trig-another}
used in the book~\cite{BatemanErdelyi},
we use it later.

We are interested in $2\pi$-periodic solutions
of the Lam\'e equation~\eqref{Lame-trig}.
Usually $0<\hat{k}<1$ and $n$ are fixed parameters
and $h$ plays the role of an eigenvalue. The
following Proposition holds.

\begin{Proposition}\label{Lame-properties}
Given $0<\hat{k}<1$ and $n,$
there exist an infinite sequence
of values
$$
h_0<h_1\leqslant h_2<h_3\leqslant h_4<\ldots
$$
of the parameter $h$ such that if $h=h_i$ then the 
the Lam\'e equation~\eqref{Lame-trig}
has a $2\pi$-periodic solution $\varphi_i(y)\ne0.$

For $h=h_0$ a solution $\varphi_0(y)$
is unique up to
a multiplication by a non-zero constant. 

If $h_{2i+1}(l)<h_{2i+2}(l),$ 
then solutions $\varphi_{2i+1}(y)$
and $\varphi_{2i+2}(y)$
are unique up to
a multiplication by a non-zero constant. 

If  $h_{2i+1}(l)=h_{2i+2}(l),$ 
then there exist two independent solutions
$\varphi_{2i+1}(y)$ and $\varphi_{2i+2}(y)$ 
corresponding to $h=h_{2i+1}=h_{2i+1}.$

The solution $\varphi_0(y)$ has no zero on $[0,2\pi).$
For $i\geqslant0$ both solutions 
$\varphi_{2i+1}(y)$ and $\varphi_{2i+2}(y)$ 
have exactly $2i+2$ zeroes on $[0,2\pi).$
\end{Proposition}

Our main interest is the case $n=1.$ In this case three
wonderful solutions of the Lam\'e equation~\eqref{Lame}
are known,
$$
\Ec_1^0(z)=\dn z,\quad\Ec_1^1(z)=\cn z,\quad\Es_1^1(z)=\sn z,
$$
where we use the notation used by Ince in the paper~\cite{Ince1}.
Using standard properties of the Jacobi elliptic functions
and change of variable~\eqref{sn-sin} we obtain
three solutions of the Lam\'e equation in the trigonometric
form~\eqref{Lame-trig},
$$
\Ec_1^0(y)=\sqrt{1-\hat{k}^2\sin^2y},\quad%
\Ec_1^1(y)=\cos y,\quad\Es_1^1(y)=\sin y.
$$

\begin{Proposition}\label{h-012}
If $n=1$ then we have
\begin{align*}
\varphi_0(y)=\Ec_1^0(y)=\sqrt{1-\hat{k}^2\cos^2y},\quad h_0=\hat{k}^2,\\
\varphi_1(y)=\Ec_1^1(y)=\cos y,\quad h_1=1,\\
\varphi_2(y)=\Es_1^1(y)=\sin y,\quad h_2=1+\hat{k}^2.
\end{align*}
\end{Proposition}

\noindent{\bf Proof.} The function $\Ec_1^0(y)=\sqrt{1-\hat{k}^2\cos^2y}$
has no zeroes, hence by Proposition~\ref{Lame-properties}
it is $\varphi_0(y).$ Direct check by substitution shows that $h_0=\hat{k}^2.$
The same argument works for $\varphi_1(y)$ and $\varphi_2(y).$
$\Box$

We should remark that in general $h_i$ are roots of a very complicated
transcendental equation with parameters $n$ and $\hat{k}$
and cannot be found explicitly.

Using the same approach as in Proposition~\ref{lambda-l}
we can prove the following Proposition.

\begin{Proposition}\label{h-k}
Let us fix $n=1$ and consider $h_3$ as a function of $\hat{k}^2,$
where $0<\hat{k}^2\leqslant1.$ Then $h_3(\hat{k}^2)$
is a decreasing function.
\end{Proposition}

When $\hat{k}=1$ the Lam\'e equation~\eqref{Lame} is
called degenerate because in this case we have $\sn z=\tanh z.$

\begin{Proposition}\label{k-1}
Let $n=1$ and $\hat{k}=1.$ Then we have
$$
h_0=h_1=1,\quad h_2=h_3=2.
$$
\end{Proposition}

\noindent{\bf Proof} follows immediately
from the explicit formulae for $h_i$
in the paper~\cite[Section~9]{Ince1}. $\Box$

Propositions~\ref{h-k} and~\ref{k-1} imply
the following Proposition.

\begin{Proposition}\label{h-3}
Let $n=1.$ Then for $0<\hat{k}^2<1$ we have
$$
h_3>2.
$$
\end{Proposition}

\section{Proof of the Theorem}

\subsection{Case of a double cover $\hat{\tau}_{m,k}$
of a Lawson torus}

It is easy to see from Proposition~\ref{1-1}
that for the  double cover $\hat{\tau}_{1,1}$ of
the Clifford torus
the statement of the Theorem holds. Hence
we can exclude this case from future considerations
and suppose that $m>k>1,$ $m\equiv k\equiv1\mod2,$ $(m,k)=1.$

The eigenvalues $\lambda_0(l)$ of the periodic
Sturm-Liouville problem~\eqref{Sturm-Liouville},
\eqref{periodic-phi} are always of multiplicity
one, see Proposition~\ref{SL-properties}.
Hence Proposition~\ref{lambda-l} implies 
that $\lambda_0(l)$ is a strictly increasing function of $l.$
Let us denote by $l_c$
the solution of the equation
$$
\lambda_0(l)=2.
$$
Then
$$
\#\{\lambda_0(l)|\lambda_0(l)<2,l>0,l\in\mathbb{Z}\}=
\lceil l_c\rceil-1,
$$
where $\lceil\cdot\rceil$ denotes the ceiling function,
i.e. $\lceil x\rceil=\min\{a\in\mathbb{Z}|a\geqslant x\}.$

Let us make now a crucial observation.
One can check by a direct
calculation that equation~\eqref{Sturm-Liouville}
could be written as the Lam\'e equation in the trigonometric
form~\eqref{Lame-trig} with
\begin{equation}\label{SL-Lame}
\hat{k}=\frac{\sqrt{m^2-k^2}}{m},\quad h=\lambda-\frac{l^2}{m^2},\quad
n(n+1)=\lambda.
\end{equation}
Let us remark that $0<\hat{k}<1$ since $m>k>1.$

It follows from~\eqref{SL-Lame} that
$\lambda=2$ corresponds to $n=1.$
Propositions~\ref{SL-properties} and~\ref{Lame-properties}
imply that $\lambda_0$ corresponds to $h_0.$ 
Hence we obtain
from~\eqref{SL-Lame} and Proposition~\ref{h-012}
the identities
\begin{equation}\label{h-l}
\frac{m^2-k^2}{m^2}=\hat{k}^2=h_0=\lambda_0(l)-%
\frac{l^2}{m^2}.
\end{equation}
We denoted the solution of the equation 
$\lambda_0(l)=2$ by $l_c$ and we obtain
from formula~\eqref{h-l} the equation
$$
\frac{m^2-k^2}{m^2}=2-\frac{l_c^2}{m^2}.
$$
It follows that
$$
l_c=m\sqrt{2-\frac{m^2-k^2}{m^2}}=\sqrt{m^2+k^2}.
$$
It is easy to see that $\sqrt{m^2+k^2}$ is not
integer because $m$ and $k$ are both odd. It follows
that 
$$
\#\{\lambda_0(l)|\lambda_0(l)<2,l>0,l\in\mathbb{Z}\}=%
\lceil l_c\rceil-1=[l_c]=%
\left[\sqrt{m^2+k^2}\,\right].
$$

Let us now consider $\lambda_1(l)$ and $\lambda_2(l).$
We proved in Proposition~\ref{lambda-2-1} that
$\cos y$ is an eigenfunction of the periodic
Sturm-Liouville problem~\eqref{Sturm-Liouville},
\eqref{periodic-phi} for $l=m$
and $\sin y$ is an eigenfunction 
of the same problem for $l=k.$
We know that $\cos y$ has 2 zeroes on $[0,2\pi).$ Hence
$\cos y$ could be either $\varphi_1(m,y)$ or
$\varphi_2(m,y).$ A similar argument shows that
either $\varphi_1(k,y)=\sin y$ or $\varphi_2(k,y)=\sin y.$

Let us suppose that $\varphi_2(m,y)=\cos y.$ 
Then $\lambda_2(m)=2$ and this eigenvalue has
multiplicity $1$ by Proposition~\ref{lambda-2-1}.
It follows from Propositions~\ref{lambda-possible}
and~\ref{lambda-l} that for $l\in(0,m]$
the function $\lambda_2(l)$ is strictly increasing.
Then $\lambda_2(k)<\lambda_2(m)=2$ and
by Proposition~\ref{SL-properties} we have the inequality
$\lambda_1(k)\leqslant\lambda_2(k)<\lambda_2(m)=2.$
This contradicts the fact that either $\lambda_1(k)=2$
or $\lambda_2(k)=2.$ The obtained contradiction
shows that $\lambda_1(m)=2$ and 
$\varphi_1(m,y)=\cos y.$ 
A~similar argument shows that $\lambda_2(k)=2$
and $\varphi_2(k,y)=\sin y.$ 

It is easy to see that these solutions
$\varphi_1(m,y)=\cos y$ and  $\varphi_2(k,y)=\sin y$
of the periodic Sturm-Liouville problem~\eqref{Sturm-Liouville},
\eqref{periodic-phi} correspond to the periodic solutions
$\varphi_1(y)=\Ec_1^1(y)=\cos y$ with $h_1=1$ and
$\varphi_2(y)=\Es_1^1(y)=\sin y$ with $h_2=1+\hat{k}^2$
of the Lam\'e equation, see Proposition~\ref{h-012}.

It follows from Propositions~\ref{lambda-possible}
and~\ref{lambda-l} that for $l\in(0,m]$
the function $\lambda_1(l)$ is strictly increasing
from $\lambda_1(0)$ to $\lambda_1(m)=2.$
This implies that
$$
\#\{\lambda_1(l)|\lambda_1(l)<2,l>0,l\in\mathbb{Z}\}=m-1.
$$
In a similar way we obtain
$$
\#\{\lambda_2(l)|\lambda_2(l)<2,l>0,l\in\mathbb{Z}\}=k-1.
$$

Let us suppose that $\lambda_3(0)\leqslant2.$ 
We know that
$\lambda_3(k)>\lambda_2(k)=2.$
It follows that there
exists some value $l_3\geqslant0$ such that $\lambda_3(l_3)=2.$
We know that $\lambda=2$ corresponds to $n=1$
and $\lambda_3$ corresponds to $h_3$ and we see
from formulae~\eqref{SL-Lame} and Proposition~\ref{h-3} that
$$
2-\frac{l_3^2}{m^2}=\lambda_3-\frac{l_3^2}{m^2}=h_3>2.
$$
This implies
$$
\frac{l_3^2}{m^2}<0,
$$
but this is impossible. Hence, $\lambda_3(0)>2.$
It follows from Proposition~\ref{SL-properties}
that $\lambda_i(l)>2$ for $i\geqslant3$ and $l\geqslant0.$

We are ready now to compute $N(2).$ Using 
formula~\eqref{N2cover} from Proposition~\ref{how-to-count}
we obtain
\begin{gather*}
N(2)=\#\{\lambda_0(0),\lambda_1(0),\lambda_2(0)\}+
2\#\{\lambda_0(l)|\lambda_0(l)<2,l>0,l\in\mathbb{Z}\}+\\
+2\#\{\lambda_1(l)|\lambda_1(l)<2,l>0,l\in\mathbb{Z}\}+%
2\#\{\lambda_2(l)|\lambda_2(l)<2,l>0,l\in\mathbb{Z}\}=\\
=3+2\left[\sqrt{m^2+k^2}\right]+2(m-1)+2(k-1)=
2\left(\left[\sqrt{m^2+k^2}\right]+m+k\right)-1.
\end{gather*}
The statement of the Theorem follows now
from Proposition~\ref{reduction}
and the following formula,
\begin{gather*}
\Lambda_{N(2)}(\hat{\tau}_{m,k})=%
\lambda_{N(2)}(\hat{\tau}_{m,k})\Area(\hat{\tau}_{m,k})=%
2\int_0^{2\pi}dx\int_{-\pi}^{\pi}p(y)dy=\\
=2\cdot 2\pi\cdot 4kE\left(i\frac{\sqrt{m^2-k^2}}{k}\right)%
=16\pi k\frac{m}{k}E\left(-\frac{\sqrt{m^2-k^2}}{m}\right)=%
16\pi mE\left(\frac{\sqrt{m^2-k^2}}{m}\right).
\end{gather*}

\subsection{Case of a Lawson torus $\tau_{m,k}$}

In this case the value of $N(2)$
follows immediately from the case of its
double cover $\hat{\tau}_{m,k}$ and formula~\eqref{N2torus} from
Proposition~\ref{how-to-count}. The area
of $\tau_{m,k}$ is just half of the area of
its double cover $\hat{\tau}_{m,k}.$
This gives the answer.

\subsection{Case of a Lawson Klein bottle $\tau_{m,k}$}

As we already discussed, we assume $m\equiv0\mod 2,$ $k\equiv1\mod2,$
$(m,k)=1.$

Let us consider the case $m>k.$ Then we have the same argument
as for the double cover $\hat{\tau}_{m,k}$ of a Lawson torus,
but we have to take into account parity of solutions.

As we know from Proposition~\ref{SL-parity}, the
solutions $\varphi_0(l,y)$ are even. This
implies that
$$
\#\{\lambda_0(l)|\lambda_0(l)<2,l>0,l\in\mathbb{Z},l\mbox{\ is even},%
\varphi_0(l,y)\mbox{\ is even}\}=%
\left\lceil\frac{l_c}{2}\right\rceil-1.
$$
Since $m$ is even, $k$ is odd, $l_c=\sqrt{m^2+k^2}$
cannot be an even integer, and we obtain
$$%
\#\{\lambda_0(l)|\lambda_0(l)<2,l>0,l\in\mathbb{Z},l\mbox{\ is even},%
\varphi_0(l,y)\mbox{\ is even}\}=%
\left[\frac{\sqrt{m^2+k^2}}{2}\right].
$$

As we know from Proposition~\ref{SL-parity}, a
solution $\varphi_i(l,y)$ is even or odd
and parity is preserved on an interval $a\leqslant l\leqslant b$
where $\lambda_i(l)$ is of multiplicity one.
It follows that $\varphi_1(l,y)$
is even for $0\leqslant l\leqslant m$ since
$\varphi_1(m,y)=\cos y$ is even and
$\varphi_2(l,y)$
is odd for $0\leqslant l\leqslant k$ since
$\varphi_2(k,y)=\sin y$ is odd.
This implies that
\begin{gather*}
\#\{\lambda_1(l)|\lambda_1(l)<2,l>0,l\in\mathbb{Z},l\mbox{\ is even},%
\varphi_1(l,y)\mbox{\ is even}\}=\\
=\{\lambda_1(2),\lambda_1(4),\dots,\lambda_1(m-2)\}=\frac{m}{2}-1
\end{gather*}
and
\begin{gather*}
\#\{\lambda_2(l)|\lambda_2(l)<2,l>0,l\in\mathbb{Z},l\mbox{\ is odd},%
\varphi_2(l,y)\mbox{\ is odd}\}=\\
=\{\lambda_1(1),\lambda_1(3),\dots,\lambda_1(k-2)\}=\frac{k-1}{2}.
\end{gather*}
Finally, using formula~\eqref{N2Klein} from
Proposition~\ref{how-to-count} we obtain
\begin{gather*}
N(2)=\#\{\lambda_0(0),\lambda_1(0)\}+\\
+2\#\{\lambda_0(l)|\lambda_0(l)<2,l>0,l\in\mathbb{Z},l\mbox{\ is even},%
\varphi_0(l,y)\mbox{\ is even}\}+\\
+2\#\{\lambda_1(l)|\lambda_1(l)<2,l>0,l\in\mathbb{Z},l\mbox{\ is even},%
\varphi_1(l,y)\mbox{\ is even}\}+\\
+2\#\{\lambda_2(l)|\lambda_2(l)<2,l>0,l\in\mathbb{Z},l\mbox{\ is odd},%
\varphi_2(l,y)\mbox{\ is odd}\}=\\
=2+2\left[\frac{\sqrt{m^2+k^2}}{2}\right]+2\left(\frac{m}{2}-1\right)+%
2\frac{k-1}{2}%
=2\left[\frac{\sqrt{m^2+k^2}}{2}\right]+m+k-1.
\end{gather*}
The statement of the Theorem follows now
from Proposition~\ref{reduction}
and the following formula,
\begin{gather*}
\Lambda_{N(2)}(\tau_{m,k})=\lambda_{N(2)}(\tau_{m,k})\Area(\tau_{m,k})=%
2\int_0^{\pi}dx\int_{-\pi}^{\pi}p(y)dy=\\
=2\cdot \pi\cdot 4kE\left(i\frac{\sqrt{m^2-k^2}}{k}\right)%
=8\pi k\frac{m}{k}E\left(-\frac{\sqrt{m^2-k^2}}{m}\right)=%
8\pi mE\left(\frac{\sqrt{m^2-k^2}}{m}\right).
\end{gather*}

The case $k>m>0$ is similar. In this case we should
use another trigonometric form of the Lam\'e equation
\begin{equation}\label{Lame-trig-another}
[1-(\hat{k}\cos y)^2]\frac{d^2\varphi}{dy^2}+%
\hat{k}^2\sin y\cos y\frac{d\varphi}{dy}+%
[h-n(n+1)(\hat{k}\cos y)^2]\varphi=0
\end{equation}
used e.g. in the book~\cite{BatemanErdelyi}.
Equation~\eqref{Lame-trig-another} could be obtained
from equation~\eqref{Lame} using the change of variable
$$
\sn z=\cos y\quad\Longleftrightarrow\quad y=\frac{\pi}{2}-\am z.
$$
Equation~\eqref{Sturm-Liouville}
could be written as the Lam\'e equation in the trigonometric
form~\eqref{Lame-trig-another} with
$$
\hat{k}=\frac{\sqrt{k^2-m^2}}{k},\quad h=\lambda-\frac{l^2}{k^2},\quad
n(n+1)=\lambda.
$$
Let us remark that $0<\hat{k}<1$ since $k>m>0.$
The rest of the proof is similar to the proof
in the case $m>k>0$ and the resulting formulae
are the same.
$\Box$

\section*{Acknowledgments}

The author is very indebted to I.~V.~Polterovich who
inspired the interest to spectral geometry.
The author is also grateful to 
H.~Volkmer for providing
the paper~\cite{Volkmer}. 

The author thanks A.P.Veselov and P.~Winternitz for fruitful
discussions.

This work was partially supported
by Russian Federation Government grant no.~2010-220-01-077, 
ag. no.~11.G34.31.0005,
by the Russian Foundation
for Basic Research (grant no.~08-01-00541
and grant no.~11-01-12067-ofi-m-2011),
by the Russian State Programme for the Support of
Leading Scientific Schools (grant no.~5413.2010.1)
and by the Simons-IUM fellowship.

\end{document}